\documentclass[conference]{IEEEtran}
\IEEEoverridecommandlockouts
\usepackage{hyperref}
\usepackage[utf8]{inputenc}
\usepackage{csquotes}
\usepackage[
backend = biber,
style=ieee, 
sorting=none, 
natbib=true, 
hyperref=true, 
date=year,
isbn=false,
maxbibnames = 3, 
dashed=false, 
]{biblatex}
\AtEveryBibitem{
    \clearfield{urlyear}
    \clearfield{note}
    \iffieldundef{pages}{}{\clearfield{doi}}
    \ifentrytype{misc}{\iflistundef{publisher}{}{\clearfield{url}}}{\clearfield{url}}
}

\usepackage{siunitx}
\newcommand{\CF}{\mathit{CF}}

\usepackage{amsmath,amssymb,amsfonts}
\usepackage{multicol}
\usepackage{multirow}
\usepackage{makecell}

\usepackage{mathtools}

\newcommand*{\tran}{^{\mkern-1.5mu\mathsf{T}}}

\usepackage{accents}
\newcommand\munderbar[1]{%
  \underaccent{\bar}{#1}}

\usepackage{graphicx}
\usepackage{tabularx}
\usepackage{textcomp}
\usepackage{xcolor}
\def\BibTeX{{\rm B\kern-.05em{\sc i\kern-.025em b}\kern-.08em
    T\kern-.1667em\lower.7ex\hbox{E}\kern-.125emX}}

\usepackage{algorithm}
\usepackage{algpseudocode}

\begin{document}

\title{Towards time series aggregation with exact error quantification for optimization of energy systems}

\author{\IEEEauthorblockN{Beltrán Castro Gómez, Yannick Werner, and Sonja Wogrin}
\IEEEauthorblockA{Institute of Electricity Economics and Energy Innovation, Graz University of Technology, Graz, Austria \\
\{castrogomez,yannick.werner,wogrin\}@tugraz.at}}

\IEEEpubid{\makebox[\columnwidth]{979-8-3315-1278-1/25/\$31.00 ©2025 IEEE\hfill}\hspace{\columnsep}\makebox[\columnwidth]{ }}

\maketitle

\IEEEpubidadjcol

\begin{abstract} 
Energy system optimization models are becoming increasingly popular for analyzing energy markets, such as the impact of new policies or interactions between energy carriers. One key challenge of these models is the trade-off between modeling accuracy and computational tractability.
A recently proposed mathematical framework addresses this challenge by achieving exact time series aggregations merging time periods sharing the same active constraint sets. This aggregation, however, is insufficient when the number of unique active constraints is large. We overcome this issue by aggregating data points from different active constraint sets. While this further reduces model size, it inevitably introduces an error compared to the full model. Yet, we show how this error can be exactly quantified without re-solving the optimization problem, enabling users to trade off computational efficiency and model accuracy proactively.
 This may be especially useful in energy markets to accommodate varying granularity across short- and long-term time horizons.
\end{abstract}


\section{Introduction}
\label{sec:section-i}
Energy system optimization models (ESOMs) represent highly complex systems in several dimensions, e.g., technological, temporal, spatial, or sectoral \cite{ridha_complexity_2020}.
The complexity has risen in recent years, e.g., due to the integration of electricity from intermittent renewable energy sources (RES) \cite{kroposki_integrating_2017}, or the increasing importance of energy markets to accommodate, e.g., an increasing number of small participants, such as residential prosumers with photovoltaic and storage systems \cite{sousa_peer--peer_2019, ilic_energy_2012}, or new and potentially nonlinear sector coupling technologies, like electrolyzers \cite{mohammadi_comprehensive_2018, robinius_power--gas_2018} or heat pumps \cite{gaur_heat_2021, abd_alla_electrification_2022}.
The main difficulty of large-scale ESOMs relies upon the trade-off between modeling detail and computational tractability \cite{brochin_harder_2024}. Complexity reduction techniques can improve computational tractability, e.g., by reducing the length of the considered time horizon. To support efficient decision-making in energy markets and systems, however, preserving the outputs of the original ESOM throughout the aggregation process is essential. While we are convinced that the fundamentals of our method extend to other complexity dimensions, we will only focus on time series aggregation (TSA) \cite{kotzur_modelers_2021} in this paper.

Among TSA methods for ESOMs, we can differentiate between two main categories: a-priori and a-posteriori methods \cite{hoffmann_review_2020}.
Traditional a-priori TSA methods try to best approximate input time series data, e.g., hourly electricity demand, by minimizing the difference between the original and the aggregated time series according to a given metric, such as Euclidean distance in k-means \cite{lloyd_least_1982}.
Such a-priori TSA methods can be improved for ESOMs, e.g., by adding \textit{extreme} periods to the aggregated time series \cite{pfenninger_dealing_2017, kotzur_impact_2018, poncelet_selecting_2017}, or including extreme values during the aggregation process itself \cite{mavrotas_mathematical_2008, green_divide_2014, kotzur_impact_2018}.
However, a-priori TSA methods may not necessarily lead to a good approximation of the ESOM outputs compared to the model with the complete time horizon -- as illustrated in \cite{wogrin_time_2023} -- because they are unaware of how the optimization model relates input data to model outputs. A good TSA method should allow one to control the input data aggregation depending on the \textit{output error}, defined between the \textit{full} and \textit{aggregated} optimization model outputs\footnote{We refer to the full model as the one that is based on the original time series, while the aggregated model uses the aggregated time series.}. Neglecting this relation may result in suboptimal solutions and poor decision-making in energy markets and systems.

On the other hand, a-posteriori TSA methods additionally consider optimization model outputs in the aggregation process. In \cite{bahl_time-series_2017}, the authors propose an iterative method that estimates investment decisions by bounding the output error of the corresponding hourly operational problem. Reference~\cite{sun_data-driven_2019} obtains representative days by clustering the cost of optimal investment decisions of individual days and \cite{hilbers_reducing_2023} samples \textit{extreme} and \textit{regular} periods, measured based on aggregated model outputs. However, we identify a gap in a-posteriori TSA: the methods that we have found in the literature conceive heuristic components in the aggregation process, e.g., \cite{sun_data-driven_2019}~relies on the assumption that similarity in objective function value relates to similarity in overall optimization outcome, and both \cite{bahl_time-series_2017, hilbers_reducing_2023}~involve preliminary a-priori clustering outcomes based solely on the input data domain. Thus, they fail to provide analytical expressions to compute exact output error bounds or provide metrics to relate the aggregation to the output error.

In previous work \cite{wogrin_time_2023}, an a-posteriori TSA framework that aggregates data points from time periods that share the same set of active constraints, referred to as \textit{basis}, was proposed\footnote{Following the basic feasible solution terminology of the Simplex method.}, and we use both terminologies interchangeably. 
Under certain assumptions for linear programs (LPs), this TSA method has been proven to be exact and yield zero error in aggregated model outputs \cite{wogrin_time_2023}. However, if the number of unique active constraint sets is large, the aggregation of data points solely within each set may only lead to a minor reduction in model size. We address this issue in this paper.  

The original contributions of this paper are: we extend the aforementioned theoretical framework with a methodology to merge data points from different active constraint sets; this allows us to derive an exact analytical expression to quantify the associated output error without having to solve the new model.
Using this expression, we elaborate on how to optimally choose the active constraint sets from which to merge data points considering the combinatorial complexity of this decision. The proposed methodology assists the modeler in efficiently navigating the trade-off between computational tractability and model accuracy.

The remainder of this paper is organized as follows: Section~\ref{sec:section-ii} recaps the TSA framework based on active constraint sets and introduces the optimal transport problem, which serves as an example. Section~\ref{sec:section-iii} extends the aforementioned TSA framework by proposing a methodology to merge data points from different active constraint sets. Section~\ref{sec:section-iv} presents and discusses our experimental results to showcase our proposed method applied to the optimal transport problem. Finally, Section~\ref{sec:section-v} concludes.

\section{Preliminaries}
\label{sec:section-ii}
First, Section~\ref{subsec:basis-def} shows how TSA can be performed based on active constraint sets. Afterward, Section~\ref{subsec:transport} introduces the optimal transport problem.

\subsection{TSA based on active constraint sets}
\label{subsec:basis-def}
Let $t \in \mathcal{T} = \{1, 2, \dots T\}$ denote the set of time steps. For the sake of simplicity, we assume that there are no time-linking constraints\footnote{In related work \cite{cardona-vasquez_enhancing_2024}, time-linking constraints such as ramping are included in the underlying data aggregation framework based on active constraint sets. Future work will address storage technologies. However, such time-linking constraints are out of the scope of this contribution.}. Let us define a standard-form primal (left) and dual (right) linear program (LP):
\vspace{-.8cm}
\begin{multicols}{2}
    \hspace{-0.35cm}
    \begin{minipage}[t]{1\columnwidth}
    \begin{subequations}
    \label{eq:lp-primal}
    \begin{align}
        \min_{x_t} \quad & \sum_{t \in \mathcal{T}} {c}\tran x_t & & \label{eq:lp-primal_a} \\
        \mathrm{s.\,t.} \quad & Ax_t \leq b_t, & & \forall t, \label{eq:lp-primal_b} \\
        & x_t \geq 0, & & \forall t, \label{eq:lp-primal_c}
    \end{align}
    \end{subequations}
    \end{minipage}
    \break
    \begin{minipage}[t]{1\columnwidth}
    \begin{subequations}
    \label{eq:lp-dual}
    \begin{align}
        \max_{y_t} \quad & \sum_{t \in \mathcal{T}} {b_t}\tran y_t & & \label{eq:lp-dual_a} \\
        \mathrm{s.\,t.} \quad & {A}\tran y_t \geq c,  && \forall t, \label{eq:lp-dual_b} \\
        & y_t \geq 0, && \forall t, \label{eq:lp-dual_c}
    \end{align}
    \end{subequations}
    \end{minipage}
\end{multicols}
\vspace{-.3cm}
where $c \in \mathbb{R}^{N}$, $x_t \in \mathbb{R}^{N}$, $A \in \mathbb{R}^{M \times N}$, $b_t \in \mathbb{R}^{M}$ and $y_t \in \mathbb{R}^{M}$. Without loss of generality, we assume that $A$ and $c$ are identical across all time steps. Moreover, we require that Problems~\eqref{eq:lp-primal} and \eqref{eq:lp-dual} are feasible and that strong duality holds. Optimal solutions to decision variables are denoted by~$\{\}^*$, e.g., $x_t \rightarrow x_t^*$.

Since Problem~\eqref{eq:lp-primal} is feasible, for every time step $t \in \mathcal{T}$, there is a set of active (also called binding) constraints $\mathcal{A}_t \subseteq \{1, 2, \dots, M\}$, that are satisfied with equality.
Then, we define $i \in \mathcal{I}$ as the set of bases, i.e., unique active constraint sets, $\mathcal{A}_i$, such that $\mathcal{A}_i \neq \mathcal{A}_{i'}$ for any $i,i' \in \mathcal{I}$ with $i \neq i'$, and $\mathcal{T}_i \subseteq \mathcal{T}$, such that $\mathcal{A}_t = \mathcal{A}_i, \quad \forall t \in \mathcal{T}_i$. 
Note that since $\mathcal{T}_i \neq \emptyset, \forall i \in \mathcal{I}$, and $\bigcap_{i\in \mathcal{I}} \mathcal{T}_i = \emptyset$, it always holds that $|\mathcal{I}| \leq |\mathcal{T}|$. However, ideally, we want $|\mathcal{I}| \ll |\mathcal{T}|$.

Using set $\mathcal{I}$ related to the unique active constraint sets, we can aggregate Problems~\eqref{eq:lp-primal} and \eqref{eq:lp-dual}, respectively,  to:
\vspace{-0.8cm}
\begin{multicols}{2}
    \hspace{-0.35cm}
    \begin{minipage}[t]{1.05\columnwidth}
    \begin{subequations}
    \label{eq:agg-lp-primal}
    \begin{align}
        \min_{\bar{x}_i} ~ & \sum_{i \in \mathcal{I}} c \tran \bar{x}_{i} W_{i} \hspace{-0.2cm} & \label{eq:agg-lp-primal_a} \\
        \mathrm{s.\,t.} ~ & A \bar{x}_{i} \leq \bar{b}_{i}, && \forall i, \label{eq:agg-lp-primal_b} \\
        & \bar{x}_i \geq 0, & & \forall i, \label{eq:agg-lp-primal_c}
    \end{align}
    \end{subequations}
    \end{minipage}
    \break
    \begin{minipage}[t]{1\columnwidth}
    \begin{subequations}
    \label{eq:agg-lp-dual}
    \begin{align}
        \max_{\bar{y}_i} ~ & \sum_{i \in \mathcal{I}} \bar{b}_{i} \tran \bar{y}_{i} \hspace{-0.2cm} & \label{eq:agg-lp-dual_a} \\
        \mathrm{s.\,t.} ~ & A\tran \bar{y}_{i} \geq c, && \forall i, \label{eq:agg-lp-dual_b} \\
        & \bar{y}_i \geq 0, && \forall i, \label{eq:agg-lp-dual_c}
    \end{align}
    \end{subequations}
    \end{minipage}
\end{multicols}
\vspace{-.3cm}
where $W_i = |\mathcal{T}_i|$ is the weight of basis $i$, i.e., the number of time steps with the same set of active constraints, and $\bar{b}_{i}~=~\frac{1}{W_i} \sum_{t \in \mathcal{T}_i}{b_t}$ the centroid of the data points in basis $i$.

As proven in \cite{wogrin_time_2023}, the optimal solutions to Problems~\eqref{eq:lp-primal} and \eqref{eq:agg-lp-primal} are identical in terms of objective function value and average of the decision variables $\frac{1}{W_i}\sum_{t \in \mathcal{T}_i} x_t^* = x_i^*, {\forall i \in \mathcal{I}}$. Considering the potentially much smaller dimension of Problem~\eqref{eq:agg-lp-primal} compared to \eqref{eq:lp-primal} when $|\mathcal{I}| \ll |\mathcal{T}|$, it is computationally preferable to solve Problem~\eqref{eq:agg-lp-primal} when averages of the optimal decision variables suffice.
Since we are assuming a constant cost function $c_t = c, \forall t \in \mathcal{T}$, all time steps $\mathcal{T}_i$ associated with active constraint set $i$ have the same optimal dual solution, such that ${y_t^* = y_i^*}, \forall t \in \mathcal{T}_i$. Therefore, when including weight $W_i$ of basis~$i$ in Equation~\eqref{eq:agg-lp-primal_a}, the dual variables of the aggregated Problem~\eqref{eq:agg-lp-dual} change such that $\bar{y}_i^* = W_i y_t^*$.

\subsection{Optimal transport problem: A case study}
\label{subsec:transport}
\begin{figure}
\vspace{-0.3cm}
\centerline{\includegraphics[width=0.8\columnwidth]{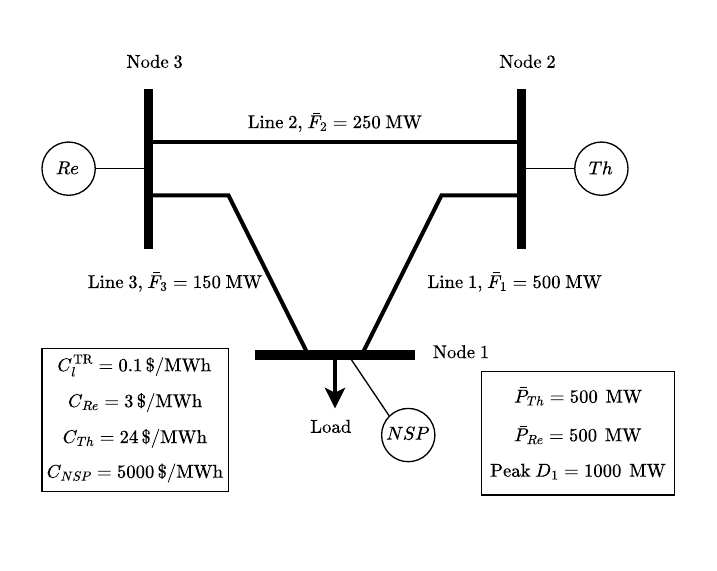}}
\caption{Stylized 3-node case study diagram including line and generator parameters for the optimal transport problem~\eqref{eq:transport-primal}.}
\label{fig:sys-diagram}
\end{figure}

Let us introduce an example of a standard LP problem in energy markets: the optimal electricity transport problem. 
We define sets of generators $g \in \mathcal{G}$, nodes $n \in \mathcal{N}$, lines $l \in \mathcal{L}$ including one line per flow direction, and time steps $t \in \mathcal{T}$, with $|\mathcal{T}| = 8736$ hours corresponding to 52 weeks.
The set $\mathcal{G}_n$ and $\mathcal{L}_n^{\mathrm{out/in}}$ collect the generators at and lines going in to (out of) node~$n$.
Let $p_{g,t}$ denote the power production of generator $g$ in time step~$t$, $f_{l,t}$ the power flow of line $l$ in time step $t$, $C_g$ the variable production cost of generator $g$, $C^{\mathrm{TR}}_{l}$ the transportation cost of line $l$, $\bar{P}_g$ the installed capacity of generator $g$, $\CF_{g,t}$ the capacity factor of generator $g$ in time step $t$, $\bar{F}_{l}$ the symmetric power flow limit of line $l$, and $D_{n,t}$ the power demand in node~$n$ in time step $t$. 
Then the optimal transport problem is given by:
\vspace{-0.1cm}
\begin{subequations}\label{eq:transport-primal}
\begin{align}
    \min_{p_{g,t}, f_{l,t}} \quad & \sum_{g,t}{C_g p_{g,t}} + \sum_{l, t}{C^{\mathrm{TR}}_{l} f_{l,t}} \hspace{-0.5cm} &&&& \label{eq:transport-primal-of} \\
    \mathrm{s.\,t.} \quad & 0 \leq p_{g,t} \leq \bar{P}_g\,\CF_{g,t} & &: \munderbar{\mu}_{g,t}, \bar{\mu}_{g,t}  & & \forall g,t, \label{eq:transport-primal-gen} \\
    & 0 \leq  f_{l,t} \leq \bar{F}_{l} & &: \munderbar{\eta}_{l,t}, \bar{\eta}_{l,t} &&  \forall l,t, \label{eq:transport-primal-flow} \\
    & D_{n,t} = \sum_{l \in \mathcal{L}^{\mathrm{in}}_n} f_{l,t} -  &&&& \notag \\ 
    & \quad \sum_{l \in \mathcal{L}^{\mathrm{out}}_n} f_{l,t} + \sum_{g \in \mathcal{G}_n}{p_{g,t}} & &: \lambda_{n, t} & & \forall n,t. \label{eq:transport-primal-balance}
\end{align}
\end{subequations}
The objective function in \eqref{eq:transport-primal-of} minimizes total system cost, which consists of generation and transmission costs. Constraints~\eqref{eq:transport-primal-gen} correspond to the lower and upper production limits of the generators. Constraints~\eqref{eq:transport-primal-flow} limit the line power flow in both directions. Constraint~\eqref{eq:transport-primal-balance} presents the nodal power balance. Based on the dual variables, defined after the colon, we derive the dual objective function of Problem~\eqref{eq:transport-primal}:
\begin{equation}
    \max_{\Theta} \quad  \sum_{g,t}{ \bar{\mu}_{g,t} \bar{P}_g \CF_{g,t} }  \,  + \sum_{l, t}{ \bar{F}_l \bar{\eta}_{l, t} } + \sum_{n,t}{\lambda_{n,t} D_{n,t}}, \label{eq:transport-dual-of}
\end{equation}
where $\Theta = \{\bar{\mu}_{g,t}, \bar{\eta}_{l, t}, \lambda_{n,t}\}$. We present the numerical results for our case study based on Equation~\eqref{eq:transport-dual-of} in Section~\ref{sec:section-iv}. We consider the 3-node system depicted in Figure~\ref{fig:sys-diagram} as a stylized but illustrative case study. In addition to generators \textit{Renewable} (\textit{Re}) and \textit{Thermal} (\textit{Th}), we also consider non-supplied power (\textit{NSP}) as a generator\footnote{We assume \textit{Th} and \textit{NSP} generators to have a capacity factor equal to 1.}.
There are two TS in the input data: the capacity factors of the wind generator, $\CF_{Re,t}$; and the demand at Node~1, $D_{1,t}$.

\begin{figure}[h]
\vspace{-0.3cm}
\centerline{\includegraphics[width=0.8\columnwidth]{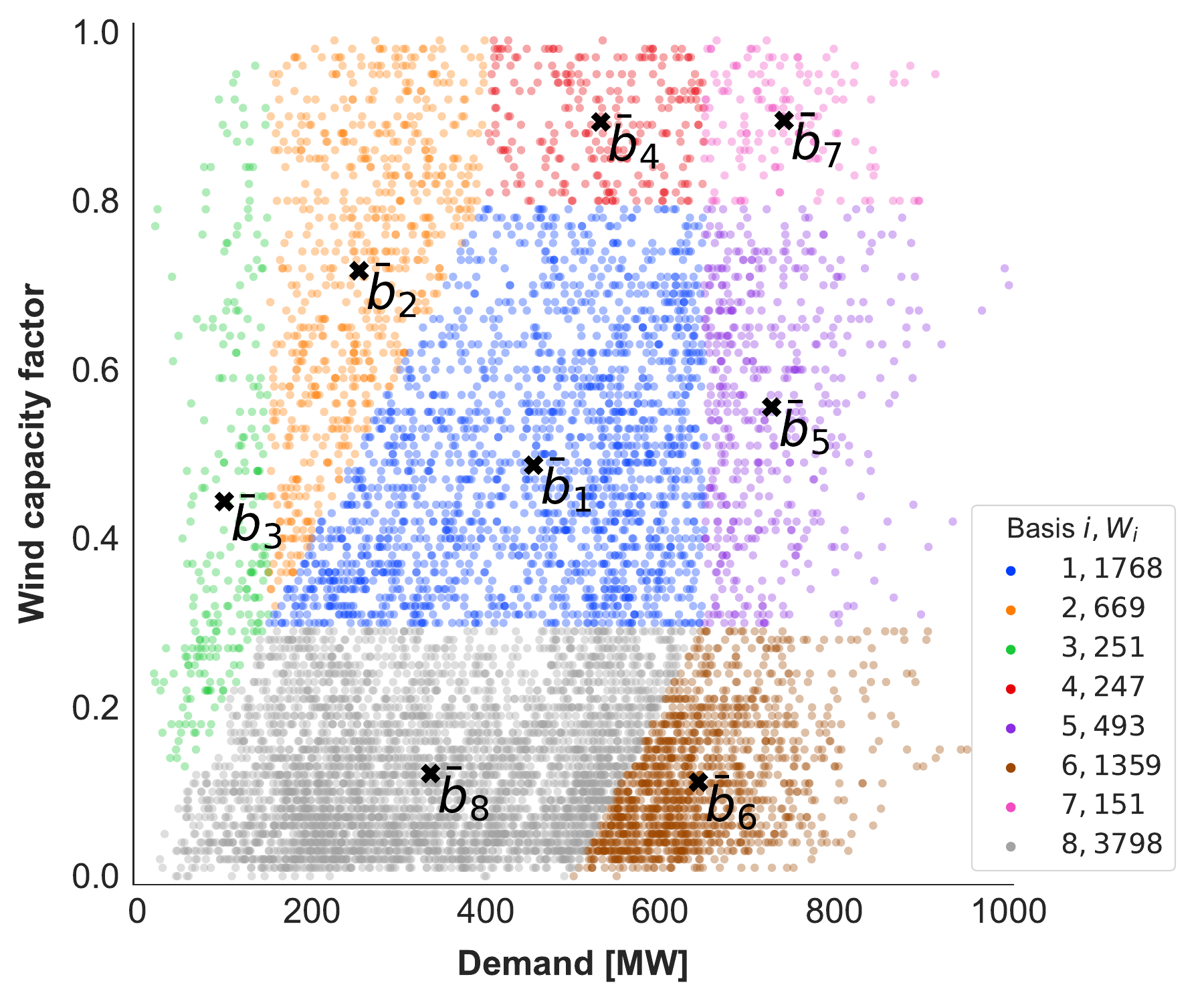}}
\caption{Hourly demand $D_{1,t}$ and wind capacity factor $\CF_{Re,t}$ time series, colored according to their basis $i$, relating input data space to model outputs. Crosses mark the centroid of each basis.}
\label{fig:ts-scatterplot-bases}
\end{figure}

After solving the optimal transport problem~\eqref{eq:transport-primal} for our case study, we observe eight unique sets of active constraints (bases), such that $\mathcal{I}=\{1,2,\dots,8\}$.
We visualize the normalized time series input data for $\CF_{Re,t}$ and $D_{1,t}$ in Figure~\ref{fig:ts-scatterplot-bases}, where every point represents a capacity factor-demand pair for a particular hour, colored according to their basis $i$.
Table~\ref{tab:bases-description} in Appendix~\ref{appendix:tables} provides intuitive descriptions for the corresponding active constraint sets. For example, for all hours in basis~2 (orange), the wind generator is the marginal generator (at node~1), and line~3 is congested.

\section{Bases Merging Framework}
\label{sec:section-iii}
In large-scale and complex modern energy systems, the cardinality of $\mathcal{I}$ may be close to that of $\mathcal{T}$ \cite{cardona-vasquez_enhancing_2024}. 
This compromises the efficiency of the aggregation and may, in the worst case, lead to an aggregated model that is still computationally intractable. Therefore, in this section, we extend the existing active constraint sets aggregation framework presented in the last section further to reduce the size of the aggregated optimization model. In particular, we describe in Section~\ref{subsec:merge-def} how we can reduce the number of bases by merging data points from different active constraint sets and, in Section~\ref{subsec:cost-misclassification}, how this impacts the approximation accuracy compared to the full model output. Note that in this paper, we assume that the mapping of time steps to bases is known through past model runs. That is the only information required to run the aggregated Problem~\eqref{eq:agg-lp-primal}. In future research, we will investigate methods to acquire the mapping without solving the full optimization problem, e.g. parametric programming.

\subsection{Definition of bases merging}
\label{subsec:merge-def}
The main idea of our proposed method is to reduce the number of bases in $\mathcal{I}$ by merging data points from different bases and replacing them with a single centroid.
With slight abuse of terminology, we refer to this as \textit{bases merging}.
In this paper, we only consider merging all data points of one basis with all data points of a single other basis. Consequently, when deciding which bases to merge, it is sufficient to only look at the centroids of the data points of each basis. While this may not necessarily be optimal, it helps us illustrate the general idea and derive some mathematical foundations of our proposed method. We will relax this assumption in future research.

\begin{figure}
\centerline{\includegraphics[width=0.8\columnwidth]{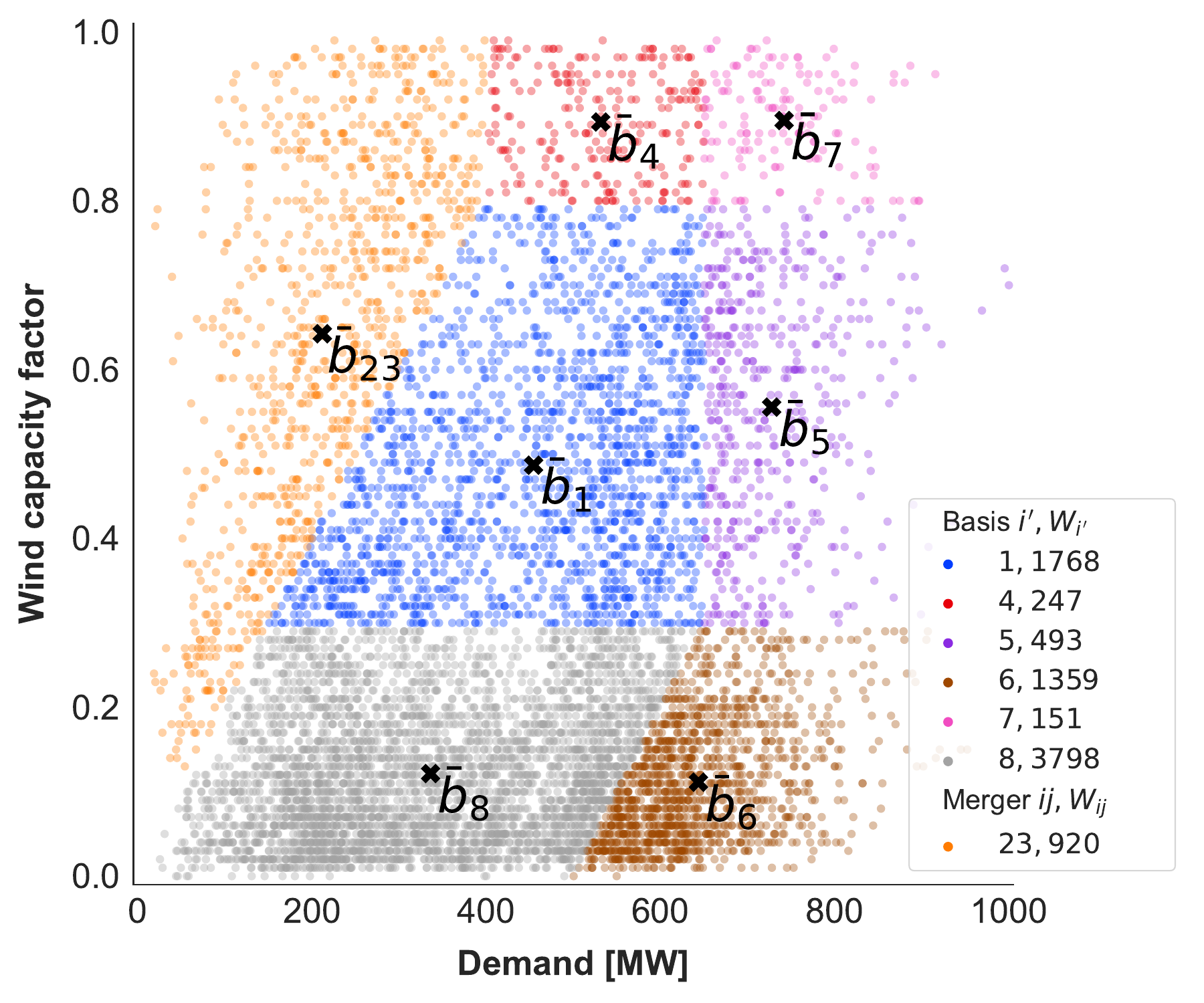}}
\vspace{-0.3cm}
\caption{Basis membership of time series tuples after merging data points from bases 2 and 3.}
\label{fig:merge-example}
\end{figure}

We define the set of \textit{bases mergers} or clusters $k \in \mathcal{K}$, such that $\mathcal{K}$ can contain elements from $\mathcal{I}$ or combinations of elements from $\mathcal{I}$, but it always holds that $|\mathcal{K}| \leq |\mathcal{I}|$.
Let us consider in the following, that we want to merge data points from, e.g., two bases $i,j \in \mathcal{I}$ with centroids $\bar{b}_i$ and $\bar{b}_j$. Then, the new centroid is given by $\bar{b}_{ij} = \frac{1}{W_i + W_j} (W_i \bar{b}_i + W_j \bar{b}_j)$.
For example, when merging all the data points from bases~2 and 3 in Figure~\ref{fig:ts-scatterplot-bases}, the new centroid will fall into basis~2.\footnote{Note that when merging data points from two bases not next to each other, e.g., 3 and 5, the new centroid may even fall outside these two bases, e.g., into basis~1. We will elaborate on this possibility in Section~\ref{subsec:cost-misclassification}.} Consequently, all data points of basis~3 will be treated as if they would fall into basis~2. This is illustrated in Figure~\ref{fig:merge-example}.
Therefore, the points initially in basis~3 will be assigned to a false set of active constraints, incurring an output error, which we discuss next.

\subsection{Cost of misclassification}
\label{subsec:cost-misclassification}
The aggregated LP associated with the merge of all points in bases $i$ and $j$ is defined as:
\vspace{-0.1cm}
\begin{subequations}
\label{eq:lp-primal-merge}
\begin{align}
    \min_{\bar{x}_{i'}, \bar{x}_{ij}} ~& \sum_{i' \in \mathcal{I'}}{c \tran \bar{x}_{i'} W_{i'}} + c \tran \bar{x}_{ij} W_{ij} & \label{eq:lp-primal-merge_a} \\
    \mathrm{s.\,t.} ~& A \bar{x}_{i'} \leq \bar{b}_{i'}, & \forall i' \in \mathcal{I'}, \label{eq:lp-primal-merge_b} \\
    & A \bar{x}_{ij} \leq \bar{b}_{ij}, & \label{eq:lp-primal-merge_c} \\
    & \bar{x}_{i'} \geq 0, & \forall i' \in \mathcal{I'}, \label{eq:lp-primal-merge_d} \\
    & \bar{x}_{ij} \geq 0, &\label{eq:lp-primal-merge_e}
\end{align}
\end{subequations}
where $\mathcal{I'} = \mathcal{I} \setminus \{i, j\}$, $\mathcal{K} = \mathcal{I}' \, \cup \{ij\}$, and $W_{ij} = W_i + W_j$ as the total weight of the new cluster, resulting from the merge of all data points from bases $i,j \in \mathcal{I}$. Since some data points in $\mathcal{T}_{ij}$ are assigned to a false set of active constraints, the optimal objective function values of Problems~\eqref{eq:agg-lp-primal} and \eqref{eq:lp-primal-merge} will not coincide. Hence, merging data points from bases $i$ and $j$ incurs an output error in the total system cost for which we derive an analytical formulation in the following. We will term this error \textit{cost of misclassification} (CoM).

Let $\mathrm{OV}^{\mathcal{I}}$ and $\mathrm{OV}^{\mathcal{K}}$ denote the objective function values of Problems~\eqref{eq:agg-lp-primal} and \eqref{eq:lp-primal-merge} at optimum, respectively.
Then, we define the cost of misclassification associated with merging points from bases $i$ and $j \in \mathcal{I}$, assuming the new centroid $\bar{b}_{ij}$ falls into basis $h \in \mathcal{I}$, as:
\vspace{-0.1cm}
\begin{equation}
\label{eq:com-ij}
    \mathrm{CoM^{ij|h}} = \mathrm{OV}^{\mathcal{I}} - \mathrm{OV}^{\mathrm{K}} = \bar{b}_{i}\tran ( \bar{y}_i^* - \frac{W_i}{W_h} \bar{y}_h^*) + \bar{b}_{j} \tran (\bar{y}_j^* - \frac{W_j}{W_h} \bar{y}_h^*).
\end{equation}
The derivation of the expression in Equation \eqref{eq:com-ij} can be found in Appendix~\ref{appendix:com}, with a generalized expression of the CoM.

We can calculate the centroid as shown, but it is not clear ex-ante in which basis the centroid will fall. In order to find out, we could, for example, check which set of active constraints is a feasible solution for the centroid. However, we have observed a favorable property of the Cost of Misclassification, which we want to mathematically prove in future research.
When computing the Cost of Misclassification for all new possible centroids of a merge, say $\mathrm{CoM^{ij | i}}$, $\mathrm{CoM^{ij | j}}$ and $\mathrm{CoM^{ij | h}}$ in our previous example, we have observed that it always holds that the new centroid of the merge coincides with the one incurring the least CoM, e.g., $\mathrm{CoM^{ij | h}} = \min_{i' \in \mathcal{I}} \mathrm{CoM}^{ij | i'}$. Therefore, we conjecture that there is no need to bound the error when we cannot know the location of the new merge, since we can know the exact error that the merge would incur.


This gives modelers control over the TSA procedure, as they can exactly quantify, without needing to solve any aggregated model -- assuming the mapping of time steps to bases is known -- the error that they are inducing in their models for any desired level of aggregation, as opposed to traditional TSA methods, that are unaware of the error that they introduce in the optimization model.

This issue also affects a-priori TSA methods. For example, when employing k-means, the obtained cluster centroids would also fall into one of the 8 bases that we have shown before. That means that all the hours that belong to that cluster, would be represented by one operational state of the power system. Just that, when doing a-priori TSA methods, one is usually ignorant of the error and misclassification one is committing.

\section{Results}
\label{sec:section-iv}

We now apply our new bases-merging methodology to the previously introduced transport problem~\eqref{eq:transport-primal} using the stylized case study in Figure~\ref{fig:sys-diagram}. In Section~\ref{subsec:enumeration}, we identify optimal bases mergers among all potential ones. Afterward, in Section~\ref{subsec:comp}, we discuss the computational complexity of finding optimal bases mergers. We use the Pyomo package v6.5.0 \cite{hart2011pyomo, bynum2021pyomo} for Python and the Gurobi Solver v11.0.0 \cite{gurobi} for our numerical analysis.
The model implementation and case study data are provided in \cite{anonymus_git}.

\subsection{Exhaustive enumeration and optimal bases mergers}
\label{subsec:enumeration}

Available bases mergers can be mathematically interpreted as partitions of a set. 
A partition of a set $\mathcal{S}$ is defined as a family of nonempty, pairwise disjoint subsets of $\mathcal{S}$ whose union is $\mathcal{S}$.\footnote{For example, the set $\{1, 2, 3\}$ has these 5 partitions: $\{ \{1\}, \{2\}, \{3\} \}$, $\{ \{1, 2\}, \{3\} \}$, $\{ \{1, 3\}, \{2\} \}$, $\{ \{1\}, \{2, 3\} \}$ and $\{ \{1, 2, 3\} \}$.}
Generally, the total number of partitions of a set of size $|\mathcal{S}| = n$ is given by the Bell Number $B_n$, which is defined recursively as $B_{n+1} = \sum_{k=0}^{n}{n \choose k} B_k$ \cite{pitman_probabilistic_1997}.

In the context of bases merging, we are looking for all partitions of the set of bases $\mathcal{I}$. In our numerical example of the transport problem, we have $|\mathcal{I}| = 8$ bases. We exhaustively enumerate all $B_8 = 4140$ possible bases mergers and solve the corresponding aggregated problems. Then, we can easily analyze how different aggregations of the input data, i.e., different bases mergers, affect the output of the aggregated model compared to the full model.
Thereby, aiming to understand good or even optimal mergers from a structural point of view.
We measure the quality of a merger in terms of the output error related to the objective function value $\varepsilon_{\mathrm{OV}}$ (with terms $\mathrm{OV}^{i}$ and $\mathrm{OV}^{k}$ defined as in Equation~\eqref{eq:com-general}, in Appendix~\ref{appendix:com}) and decision variables $\varepsilon_{g}, \forall g \in \mathcal{G}$ as:
\begin{align}
    \varepsilon_{\mathrm{OV}} & = \frac{ \sum_{ i \in \mathcal{I} }{ \mathrm{OV^i} } - \sum_{ k \in \mathcal{K} }{ \mathrm{ OV^{k} } } }{ \sum_{ i \in \mathcal{I} }{ \mathrm{OV^i} } }, \label{eq:error-of} \\
    \varepsilon_{g} &= \frac{ \sum_{ i \in \mathcal{I} }{ p_{g,i} } - \sum_{ k \in \mathcal{K} }{ p_{g,k} } }{ \sum_{ i \in \mathcal{I} }{ p_{g,i} } }. \label{eq:error-g}
\end{align}
\vspace{-0.35cm}
\begin{table}
\renewcommand{\arraystretch}{1.2} 
\renewcommand{\tabcolsep}{2.5pt}
\caption{Results of optimal bases merging. Relative errors $\varepsilon$ are shown in~$\%$. Values are rounded to two decimals.}
\centering
\begin{tabular}{|c|c|c|c|c|c|}
\hline
$|\mathcal{K}|$ & Optimal bases mergers &  $\varepsilon_{\mathrm{OV}}$ & $\varepsilon_{Re}$ & $\varepsilon_{Th}$ & $\varepsilon_{NSP}$ \\
\hline
8 & \{1\},\{2\},\{3\},\{4\},\{5\},\{6\},\{7\},\{8\} & $0.00$ & $0.00$ & $0.00$ & $0.00$ \\
\hline
7 & \{1\},\{2,\,3\},\{4\},\{5\},\{6\},\{7\},\{8\} & $ 0.00$ & $0.00$ & $0.00$ & $0.00$ \\
\hline
6 &  \{1,\,8\},\{2,\,3\},\{4\},\{5\},\{6\},\{7\} & $0.00$ & $0.00$ & $0.00$ & $0.00$ \\
\hline
5 & \{1,\,8\},\{2,\,3\},\{4\},\{5,\,7\},\{6\} & $0.02$ & $-0.59$ & $0.31$ & $0.00$ \\
\hline
4 & \{1,\,4,\,8\},\{2,\,3\},\{5,\,7\},\{6\} & $0.05$ & $-1.52$ & $0.80$ & $0.00$ \\
\hline
3 & \{1,\,2,\,3,\,4,\,8\},\{5,\,7\},\{6\} &  $0.28$ & $-9.74$ & $5.08$ & $0.00$ \\
\hline
2 & \{1,\,2,\,3,\,4,\,5,\,7,\,8\},\{6\} &  $28.53$ & $-9.74$ & $2.86$ & $30.39$ \\
\hline
1 & \{1,\,2,\,3,\,4,\,5,\,6,\,7,\,8\} &  $93.24$ & $-9.74$ & $-2.22$ & $100$ \\
\hline
\end{tabular}
\label{tab:merge-results}
\end{table}

The optimal bases mergers and the associated output errors\footnote{The mergers in Table~\ref{tab:merge-results} are optimal in terms of the error in the objective function value $\varepsilon_{\mathrm{OV}}$. Depending on the model's purpose, the metric of interest may be a different variable. For example, if the modeler is a wind power investor, the metric of interest may be the error in wind production $\varepsilon_{Re}$, while for a grid operator, it might be error in non-supplied power $\varepsilon_{NSP}$. In summary, the modeler's point of view determines the optimal mergers.} with respect to the full hourly model are shown in Table~\ref{tab:merge-results}.
We observe when reducing the number of clusters, e.g., from 8 to 7, it is optimal to merge data points from bases~ 2 and 3. In bases 2 and 3, the wind generator is the marginal unit. But they differ in the congestion of line 3, resulting in a difference in their locational marginal prices of 0.1 \$/MWh due to the transport cost $C^{\mathrm{TR}}_{l}$ (see Table \ref{tab:bases-description} in Appendix~\ref{appendix:tables}). As the transport cost is comparably small, the error in objective function value induced by the aggregation is $\varepsilon_{\mathrm{OV}} \approx 0$. 
Notably, in our example, it is possible to reduce the number of clusters down to four, without incurring an absolute relative error of above $0.05\,\%$ in the optimal objective function value and $2\,\%$ in the decision variable values. This shows the potential of the proposed bases merging framework to choose input data aggregations that reduce model complexity without losing a lot of accuracy. While this accuracy can potentially also achieved with a-priori methods, there is no way of controlling the output error during the aggregation.

Note that it is not needed for quantifying the error, to actually re-solve the aggregated optimization problems after merging bases, as the CoM \eqref{eq:com-general} given in Appendix~\ref{appendix:com}, which can be computed as shown in Section~\ref{subsec:cost-misclassification}, are equal to $\varepsilon_{\mathrm{OV}}$ in absolute terms.
However, when the number of bases $\mathcal{I}$ is large this may still be computationally challenging. We tackle this issue in the next Section.

\subsection{Computational complexity of finding optimal bases mergers}
\label{subsec:comp}
The Bell Number is growing faster than exponentially with the size of the set, in our case the bases $\mathcal{I}$, to be partitioned.
While for $n = 8$ there are 4140 partitions, for $n = 15$ it would already be 1.3 Billion. In our numerical example, however, we observed that the optimal bases mergers share some common properties drastically reducing the number of partitions to consider.
We explain both properties in the following. Due to the scope of this paper, we must postpone mathematical proofs for their validity to future research. 

\paragraph{Greediness}
\label{para:greediness}
Looking at Table~\ref{tab:merge-results}, it can be noted that for any number of bases in the merged set, the optimal bases mergers follow a hierarchical pattern. For example, merging data points from bases~2 and 3 is not only optimal when reducing the number of clusters by 1, but also for all further reductions. We therefore hypothesize that a greedy strategy will generally deliver good, if not optimal, bases mergers for all types of LPs.

\paragraph{Adjacency}
\label{para:adjacency}
We have further observed that all optimal bases mergers shown in Table~\ref{tab:merge-results} include bases that are \textit{adjacent}, meaning they share only one edge in the input space shown in Figure~\ref{fig:ts-scatterplot-bases}. For instance, bases 2 and 3 are adjacent as opposed to bases 2 and 6. 
We hypothesize that adjacency in input space relates to the adjacency of basic solutions in the Simplex method. This may be the foundation for proving it as a necessary condition for optimal bases mergers.

Table~\ref{tab:strategy-comparison-non-cummulative} in Appendix~\ref{appendix:strategies} shows the number of potential optimal bases mergers based on the exhaustive enumeration baseline, a greedy strategy, and the latter with an adjacency requirement.

Our results suggest both properties, greediness and adjacency, significantly reduce potential optimal bases mergers compared to the Bell Number, possibly overcoming its faster-than-exponential growth and achieving polynomial complexity.
In a non-degenerate LP with $m$ constraints and $n$ variables, the number of basic feasible solutions is at most $n \choose m$, which grows as $\mathcal{O}(2^n)$ in the worst case. On the other hand, each basic feasible solution has at most $m(n-m)$ adjacent basic feasible solutions, assuming there is no degeneracy and all pivots lead to feasible solutions. This grows as $\mathcal{O}(n^2)$ in the worst case. Therefore, we can expect greater reductions in larger problems as the difference of both terms grows exponentially.

\section{Conclusions}
\label{sec:section-v}
Recently, it was proposed in~\cite{wogrin_time_2023} to use TSA based on active constraint sets (or bases). This aggregation becomes inefficient when the number of active constraint sets becomes large.  We extend this framework by introducing the concept of merging data belonging to different active constraint sets, which we term bases merging. While this reduces model size it incurs an output error.
However, we demonstrate that optimal mergers may be systematically identified and derive an exact analytical formulation to compute the corresponding output error that does not require re-solving the optimization model. This enables modelers to balance model accuracy and computational complexity efficiently.

Future research should derive mathematical proofs or guarantees for choosing optimal bases mergers using a greedy strategy or requiring adjacency. Furthermore, the proposed bases merging method should be extended to account for time-linking constraints. Finally, it should be investigated how the assignment of data points to bases can be derived or approximated without solving the full model.

\section*{Acknowledgment}
Funded by the European Union (ERC, NetZero-Opt, 101116212). Views and opinions expressed are however those of the authors only and do not necessarily reflect those of the European Union or the European Research Council. Neither the European Union nor the granting authority can be held responsible for them.

\printbibliography

\clearpage

\appendices

\section{Bases of the optimal transport problem}\label{appendix:tables}
In Table~\ref{tab:bases-description}, we intuitively describe the eight bases of the optimal transport problem applied to the case study in Section~\ref{subsec:transport}. For every basis $i \in \mathcal{I}$, we present its color used in Figure~\ref{fig:ts-scatterplot-bases}, the number of data points that share this basis (weight), which network lines are congested, the generators operating at full load, and, lastly, the marginal generator (i.e., most expensive generator needed to balance supply and demand) and locational marginal price (LMP) at Node~1, where all the demand is.

\begin{table}[htbp]
\renewcommand{\arraystretch}{1.2} 
\renewcommand{\tabcolsep}{2.5pt}
\caption{Description of set of bases $\mathcal{I}$ from transport problem.}
\centering
\begin{tabular}{|c|c|c|c|c|c|c|c|}
\hline
Basis $i$ & Color & $|\mathcal{T}_i|$ & Congestion & Full load & Marginal gen. & LMP \\
\hline
1 & Blue & 1768 & Line 3 & \textit{Re} & \textit{Th} & 24.1 \\
\hline
2 & Orange & 669 & Line 3 & - & \textit{Re} & 3.2 \\
\hline
3 & Green & 251 & - & - & \textit{Re} & 3.1 \\
\hline
4 & Red & 247 & Lines 2, 3 & - & \textit{Th} & 24.1 \\
\hline
5 & Purple & 493 & Lines 1, 3 & \textit{Re} & \textit{NSP} &  5000 \\
\hline
6 & Brown & 1359 & Line 1 & \textit{Re}, \textit{Th} & \textit{NSP} & 5000 \\
\hline
7 & Pink & 151 & Lines 1, 2, 3 & - & \textit{NSP} & 5000 \\
\hline
8 & Gray & 3798 & - & \textit{Re} & \textit{Th} & 24.1 \\
\hline
\end{tabular}
\label{tab:bases-description}
\end{table}

\section{Cost of Misclassification}
\label{appendix:com}
Let $\mathrm{OV^{ij|h}}$ in Equation~\eqref{eq:ov-ij} denote the part of the dual objective function related to Constraint~\eqref{eq:lp-primal-merge_c} of Problem~\eqref{eq:lp-primal-merge} at optimum.
Since we are assuming that the centroid $\bar{b}_{ij}$ falls into basis $h \in \mathcal{I}$, the associated dual variables $\bar{y}_{ij|h}^*$ correspond to $\bar{y}_h^*$, associated to basis $h$, re-weighted according to the weight $W_{ij}$ of cluster $ij$.
\begin{subequations}
\label{eq:ov-ij}
\begin{align}
    \mathrm{OV^{ij|h}} & = \bar{b}_{ij} \tran \bar{y}_{ij|h}^* = \label{eq:ov-ij_a} \\
    & = (\frac{W_i}{W_i + W_j} \bar{b}_i \tran + \frac{W_j}{W_i + W_j} \bar{b}_j \tran) (W_i+W_j)\frac{\bar{y}_h^*}{W_h} = \label{eq:ov-ij_b} \\
    & = (\frac{W_i}{W_h} \bar{b}_i \tran + \frac{W_j}{W_h} \bar{b}_j \tran) \bar{y}_h^*. \label{eq:ov-ij_c}
\end{align}
\end{subequations}
Further, define $\mathrm{OV}^{i'} = c\tran \bar{x}_{i'}^* W_{i'} = \bar{b}_{i'}\tran \bar{y}_{i'}^* $ as the part of the primal and dual objective functions~\eqref{eq:agg-lp-primal_a} and \eqref{eq:agg-lp-dual_a}, respectively, attributed to basis $\forall i' \in \mathcal{I}$. We derive the expression of the CoM from Equation~\eqref{eq:com-ij} as: 
\begin{subequations}
\label{eq:com-ij-ext}
\begin{align}
    \mathrm{CoM^{ij|h}} & = \mathrm{OV}^{\mathcal{I}} - \mathrm{OV}^{\mathrm{K}} = \label{eq:com-ij-ext_a} \\ 
    & = \mathrm{OV^i} + \mathrm{OV^j} - \mathrm{OV^{ij|h}} = \label{eq:com-ij-ext_b} \\
    & = \bar{b}_i \tran \bar{y}_i^* + \bar{b}_j \tran \bar{y}_j^* - (\frac{W_i}{W_h} \bar{b}_i \tran + \frac{W_j}{W_h} \bar{b}_j \tran) \bar{y}_h^* = \label{eq:com-ij-ext_c} \\
    & = \bar{b}_{i}\tran ( \bar{y}_i^* - \frac{W_i}{W_h} \bar{y}_h^*) + \bar{b}_{j} \tran (\bar{y}_j^* - \frac{W_j}{W_h} \bar{y}_h^*). \label{eq:com-ij-ext_d}
\end{align}
\end{subequations}
In order to define a generalized expression of the CoM, let us consider the primal and dual LPs associated to the aggregation of data points with clusters $k \in \mathcal{K}$ as:
\begin{multicols}{2}
    \hspace{-0.35cm}
    \begin{minipage}[t]{1.05\columnwidth}
    \begin{subequations}
    \label{eq:lp-merge-general-primal}
    \begin{align}
        \min_{\bar{x}_k} ~ & \sum_{k \in \mathcal{K}} c \tran \bar{x}_{k} W_{k} \hspace{-0.2cm} & \label{eq:lp-merge-general-primal_a} \\
        \mathrm{s.\,t.} ~ & A \bar{x}_{k} \leq \bar{b}_{k}, && \forall k, \label{eq:lp-merge-general-primal_b} \\
        & \bar{x}_k \geq 0, & & \forall k, \label{eq:lp-merge-general-primal_c}
    \end{align}
    \end{subequations}
    \end{minipage}
    \break
    \begin{minipage}[t]{1\columnwidth}
    \begin{subequations}
    \label{eq:lp-merge-general-dual}
    \begin{align}
        \max_{\bar{y}_k} ~ & \sum_{k \in \mathcal{K}} \bar{b}_{k} \tran \bar{y}_{k} \hspace{-0.2cm} & \label{eq:lp-merge-general-dual_a} \\
        \mathrm{s.\,t.} ~ & A\tran \bar{y}_{k} \geq c, \hspace{-0.1cm} && \forall k, \label{eq:lp-merge-general-dual_b} \\
        & \bar{y}_k \geq 0, && \forall k. \label{eq:lp-merge-general-dual_c}
    \end{align}
    \end{subequations}
    \end{minipage}
\end{multicols}
Note that, in the case of perfect aggregation, $\mathcal{K} = \{\{i\} | i \in \mathcal{I}\}$, where every cluster contains points from only one basis. In the case from LP~\eqref{eq:lp-primal-merge}, $\mathcal{K} = \{\{ij\}\} \cup \{\{i'\} | i' \in \mathcal{I} \setminus \{i,j\}\}$, thus points from all bases except for $i$ and $j$ are represented by their original basis' centroid. It can be noted that LP~\eqref{eq:agg-lp-primal} and LP~\eqref{eq:lp-merge-general-primal}, as well as LP~\eqref{eq:agg-lp-dual} and LP~\eqref{eq:lp-merge-general-dual}, are equivalent from an optimization point of view. They differ in the way the input data centroids $\bar{b}_i$ and $\bar{b}_k$ are computed.

Now, let $\mathrm{OV}^{\mathcal{I}}$ and $\mathrm{OV}^{\mathcal{K}}$ denote the objective function values of Problems~\eqref{eq:agg-lp-primal} and \eqref{eq:lp-merge-general-primal} at optimum, respectively.
Define $\mathrm{OV}^{i} = c\tran \bar{x}_{i}^* W_{i} = \bar{b}_{i}\tran \bar{y}_{i}^* $ as the part of the primal and dual objective functions~\eqref{eq:agg-lp-primal_a} and \eqref{eq:agg-lp-dual_a}, respectively, attributed to basis $\forall i \in \mathcal{I}$.
Additionally, define $\mathrm{OV}^{k} = c\tran \bar{x}_{k}^* W_{k} = \bar{b}_{k}\tran \bar{y}_{k}^* $ as the part of the primal and dual objective functions~\eqref{eq:lp-merge-general-primal_a} and \eqref{eq:lp-merge-general-dual_a}, respectively, attributed to cluster $\forall k \in \mathcal{K}$.
Then, let us define the CoM for a general case where we merge bases from $i \in \mathcal{I}$ into clusters $k \in \mathcal{K}$ as the difference between $\mathrm{OV}^{\mathcal{I}}$ and $\mathrm{OV}^{\mathcal{K}}$:
\begin{subequations}
\label{eq:com-general}
\begin{align}
    \mathrm{CoM} & = \mathrm{OV^{\mathcal{I}}} - \mathrm{OV^{\mathcal{K}}} = \label{eq:com-general_a} \\
    & = \sum_{ i \in \mathcal{I} }{ \mathrm{OV^i} } - \sum_{ k \in \mathcal{K} }{ \mathrm{OV^k} } = \label{eq:com-general_b} \\
    & = \sum_{i \in \mathcal{I}} \bar{b}_{i} \tran \bar{y}_{i}^* - \sum_{k \in \mathcal{K}} \bar{b}_{k} \tran \bar{y}_{k}^*. \label{eq:com-general_c}
\end{align}
\end{subequations}

\section{Bases merging strategies}
\label{appendix:strategies}
In Table~\ref{tab:strategy-comparison-non-cummulative}, we observe the number of potential bases mergers that are evaluated by three merging strategies (Exhaustive, Greedy, Greedy \& Adjacent) as a function of the cardinality of the set of clusters $\mathcal{K}$. Note that $i,i'$ denote bases $\in \mathcal{I}$, and $k,k' \in \mathcal{K}$ denote clusters. A cluster can contain only one basis or more than one merged bases.

All three strategies initialize $\mathcal{K}$ by assigning each of the eight bases from $\mathcal{I}$ into its own cluster. Then, they iteratively reduce the cardinality of $\mathcal{K}$ by 1, as they attempt to find at each iteration the bases merger that minimizes a given metric, i.e., the CoM in our case. For the case study, all three strategies found the same optimal bases mergers presented in Table~\ref{tab:merge-results}.

\begin{table}
\renewcommand{\arraystretch}{1.2}
\caption{Number of potential mergers for different strategies as a function of the cardinality of $\mathcal{K}$.}
\centering
\begin{tabular}{|c|c|c|c|c|c|c|c|c|}
\cline{2-9}
\multicolumn{1}{c|}{} & \multicolumn{8}{c|}{$|\mathcal{K}|$} \\
\cline{2-9}
\multicolumn{1}{c|}{} & 8 & 7 & 6 & 5 & 4 & 3 & 2 & 1 \\
\hline
Exhaustive & 1 & 28 & 266 & 1050 & 1701 & 966 & 127 & 1 \\
\hline
Greedy (G.) & 1 & 28 & 21 & 15 & 10 & 6 & 3 & 1 \\
\hline
G. \& Adjacent & 1 & 11 & 10 & 8 & 7 & 4 & 3 & 1 \\
\hline
\end{tabular}
\label{tab:strategy-comparison-non-cummulative}
\end{table}
The Exhaustive strategy is presented in Algorithm~\ref{alg:exhaustive}. At each iteration $iter$, the number of bases mergers evaluated by the strategy is given by the Stirling number of the second kind $S_{|\mathcal{I}|,iter}$, which describe the number of partitions of size $iter$ of the set of bases $\mathcal{I}$, as $S_{|\mathcal{I}|,iter} = \frac{1}{iter!} \sum_{n=0}^{iter}{ {iter \choose n} (-1)^{iter-n} n^{|\mathcal{I}|} } $ \cite{Boyadzhiev01102012}.

In this strategy, the iterations can be interpreted as independent from each other, in the sense that the merger chosen at a given iteration does not affect the merger to be chosen in the subsequent ones, as all potential mergers are evaluated.

Since a priori, we do not know what is the optimal number of clusters for a particular problem, we try out all options as calculated in Algorithm~\ref{alg:exhaustive}, which yields the Bell number, e.g. $B_8 = \sum_{iter=1}^8 S_{8,iter}$ in our case.

This strategy serves as a baseline to the Greedy and Greedy \& Adjacent strategies, as it provides an upper bound of the number of potential bases mergers.

\begin{algorithm}
\caption{Exhaustive strategy}\label{alg:exhaustive}
\begin{algorithmic}
\For{$iter = |\mathcal{I}|, |\mathcal{I}| - 1, \ldots, 1$}
    \State $mergers \gets partitions(\mathcal{I}, size=iter)$
    \State Choose merger $\in mergers$ with min. CoM
\EndFor
\end{algorithmic}
\end{algorithm}

The Greedy strategy, presented in Algorithm~\ref{alg:greedy}, can be seen as an improvement with respect to the Exhaustive strategy. In this strategy, we keep a memory of the optimal merger found in the previous step. Therefore, the number of bases mergers to evaluate at iteration $iter$ depends on the decision previously taken at $iter+1$.

At each iteration, Algorithm~\ref{alg:greedy} evaluates all possible mergers given by the pairwise combinations of the elements in $\mathcal{K}$. The merger that minimizes the CoM, we refer to it as $\{k k'\}$, is chosen and $\mathcal{K}$ is updated with this choice. This means that the original separate clusters $\{k\}$ and $\{k'\}$ are removed from $\mathcal{K}$ and replace by the merged cluster $\{k k'\}$. Once two clusters have been merged, this merge cannot be undone in subsequent iterations.

The number of bases mergers evaluated by this strategy at iteration $iter$ is equal to $ {iter + 1 \choose 2} \textrm{ if } iter < |\mathcal{I}| \textrm{ else } 1$. 

\begin{algorithm}
\caption{Greedy strategy}\label{alg:greedy}
\begin{algorithmic}
\State $\mathcal{K} \gets \{\{i\} | i \in \mathcal{I}\}$
\For{$iter = |\mathcal{I}|, |\mathcal{I}| - 1, \ldots, 1$}
    \State $mergers \gets pairs(\mathcal{K})$
    \State Choose merger $\{kk'\} \in mergers$  with min. CoM 
    \State Update $\mathcal{K} \gets \mathcal{K} \backslash \{ \{ k \}, \{k'\} \}  \cup \{k k'\} $    
\EndFor
\end{algorithmic}
\end{algorithm}
The Greedy \& Adjacent strategy, described in Algorithm~\ref{alg:greedy-adj}, can be understood as an improvement with respect to the Greedy strategy. Here, on top of the memory of optimal mergers from previous iterations, we only evaluate mergers that satisfy an adjacency condition. In the case study, we have identified the following pairs of bases as adjacent: \textit{adj} = [(1, 2), (1, 4), (1, 5), (1, 8), (2, 3), (2, 4), (3, 8), (4, 7), (5, 6), (5, 7), (6, 8)] from Figure~\ref{fig:ts-scatterplot-bases}.

For a pair of clusters $k,k' \in \mathcal{K}$ to be considered adjacent:
\begin{itemize}
    \item If each of the two clusters coincides with a basis (i.e. $k=i$ and $k'=i'$), they are adjacent if they are one of the identified pairs of adjacent bases in $adj$.
    \item If any of the two clusters $k,k'$ that are considered for a merge, contain more than one basis (e.g. $k=\{i i'\}$), they are adjacent if the resulting merged cluster can be obtained from the union of pairs of adjacent bases.
    As an example, clusters $\{2,3\}$ and $\{4\}$ are adjacent since bases 2 and 3, as well as 3 and 4, are adjacent. As opposed to, for example, clusters $\{2,3\}$ and $\{6\}$, which do not satisfy the adjacency condition since basis 6 is neither adjacent to bases 2 or 3.
\end{itemize}

Since adjacency is problem-dependent, it is not trivial to derive a closed-form expression of the number of mergers that are evaluated by this strategy. However, the Greedy strategy serves as an upper bound to it.

\begin{algorithm}
\caption{Greedy \& Adjacent strategy}\label{alg:greedy-adj}
\begin{algorithmic}
\State $\mathcal{K} \gets \{\{i\} | i \in \mathcal{I}\}$
\State $adj \gets \textrm{List of adjacent bases}$
\For{$iter = |\mathcal{I}|, |\mathcal{I}| - 1, \ldots, 1$}
    \State $mergers \gets pairs(\mathcal{K})$
    \State $adjmergers \gets filter(mergers, adj)$
    \State Choose merger $\{kk'\} \in adjmergers$ with min. CoM
    \State Update $\mathcal{K} \gets \mathcal{K} \backslash \{ \{ k \}, \{k'\} \}  \cup \{k k'\} $  
\EndFor
\end{algorithmic}
\end{algorithm}

\end{document}